\documentclass[USenglish,autoref,numberwithinsect]{lipics-v2021}
\nolinenumbers
\hideLIPIcs

\usepackage{stmaryrd} 

\newcommand{\Om}{{\boldsymbol{\omega}}}
\newcommand{\N}{\mathbb{N}}
\newcommand{\Q}{\mathbb{Q}}
\renewcommand{\P}{\mathcal{P}}
\newcommand{\DT}{\mathcal{D}_{\operatorname{T}}}
\newcommand{\DANM}{\mathcal{D}_{\operatorname{ANM}}}
\newcommand{\op}{\operatorname{op}}
\renewcommand{\implies}{\Rightarrow} 
\renewcommand{\iff}{\Leftrightarrow} 
\newcommand{\rzimplies}{\Rrightarrow}
\newcommand{\inject}{\hookrightarrow}
\newcommand{\id}{\operatorname{id}}
\newcommand{\coT}{\operatorname{coT}}

\newcommand{\Error}{\operatorname{Error}}
\newcommand{\Psymb}{\mathtt{P}}
\newcommand{\Qsymb}{\mathtt{Q}}
\newcommand{\Class}{\Omega_{\lnot\lnot}} 

\title{An order-reversing embedding of Turing degrees into Arthur–Nimue–Merlin degrees}

\author{Jean {Abou Samra}}{Eötvös Loránd University (ELTE), Budapest, Hungary \and \url{https://jean.abou-samra.fr}}{jean@abou-samra.fr}{https://orcid.org/0009-0006-1438-2427}{Funded by the European Union (ERC, HOTT, 101170308). Views and opinions expressed are however those of the author(s) only and do not necessarily reflect those of the European Union or the European Research Council. Neither the European Union nor the granting authority can be held responsible for them.}
\author{David Alexander Madore}{Télécom Paris (LTCI) and Institut Polytechnique de Paris, Palaiseau, France \and \url{http://www.madore.org/~david/}}{david.madore@telecom-paris.fr}{}{}

\authorrunning{J. Abou Samra and D. A. Madore}

\Copyright{Jean Abou Samra and David A. Madore}



\relatedversiondetails{Full version}{https://arxiv.org/abs/2603.19946}
\acknowledgements{We thank Andrej Bauer, Laurent Bienvenu, Joel David Hamkins, James E. Hanson, Takayuki Kihara, Ulysse Léchine, Mabel Najdovski, the pseudonymous MathOverflow user “Pastebee”, Arno Pauly, Andrei Romashchenko, Will Sawin, Sasha Shen, and anonymous reviewers for helpful discussions and comments. No AI tool was used in this research or in the writing of this paper.}
\begin{document}

\maketitle

\begin{abstract}
The Arthur–Nimue–Merlin degrees are a generalization of the Turing degrees introduced by Kihara as a tangible description of the partially ordered set of Lawvere–Tierney topologies on the effective topos (equivalently, subtoposes of the effective topos). They are defined in terms of a three-player game that introduces both angelic and demonic non-determinism into oracle queries. We construct an order embedding of the Turing degrees with their order reversed into the Arthur–Nimue–Merlin degrees, whose image we call the “co-Turing degrees”; we then study the order relationship of these co-Turing degrees with the (naturally embedded) Turing degrees within the Arthur–Nimue–Merlin degrees.
\end{abstract}

\section{Introduction}

We begin with a riddle. King Arthur is desperate to know whether the Grail is in a particular French castle. He can ask arbitrary questions to the all-knowing mage Merlin, who promises to answer truthfully, but with a twist. Instead of merely answering “yes” or “no” in plain English, the mischievous Merlin provides a Turing machine, which halts if and only if the answer is “yes”. Can Arthur, a mere mortal limited by the Church–Turing thesis, find the Grail's location before the end of days?

Indeed he can, with a simple trick: first ask Merlin whether the Grail is in the castle, then ask whether it is \emph{not} in the castle. Exactly one of the two resulting Turing machines halts, and he can run them in parallel to know which one, revealing the answer.

Let us spice up the rules. Each of Merlin's replies now consists of \emph{two} Turing machines; for a “yes” answer, both or none halt, whereas for a “no” answer, exactly one halts.\footnote{We are indebted to Will Sawin for suggesting this particularly elegant variant of the riddle.} It turns out that, under these rules, no strategy can ensure that Arthur learns the answer to an arbitrary boolean question in finite time. A proof of this fact (which we give in \autoref{remark-answer-to-riddle}) seems to require the main technical result of this paper (\autoref{main-strategy-result}), which shows that asking the question then its negation is essentially the only useful trick in a precise sense.

\subparagraph*{Context}
Our object of study is an extension of Turing degrees where the oracles can feature certain forms of non-determinism, as in this riddle. We call these generalized degrees the \emph{Arthur–Nimue–Merlin degrees}, owing to their explicit description by Kihara using a three-player game whose players are conventionally named Arthur, Nimue and Merlin \cite{KiharaLT}.\footnote{In \cite{KiharaOracles}, Kihara uses the term “extended Turing–Weihrauch degrees”. We prefer to eschew the association with Weihrauch's name since we never deal with oracles that can only be queried once.} Whereas Turing degrees are defined as a quotient of oracles which are functions $\N \to \N$, Arthur–Nimue–Merlin degrees are a quotient of a richer set of oracles, namely functions $\N \to \P(\P(\N))$. Informally, when a program that uses $g : \N \to \P(\P(\N))$ as an oracle makes a query $n \in \N$ to the oracle, a set $P \in g(n)$ is chosen in a benevolent way (“angelic non-determinism”), then a value $v \in P$ is chosen in an adversarial way (“demonic non-determinism”), and the program receives the answer $v$. The mathematical meaning to “benevolent” and “adversarial” is made precise by introducing a game where Arthur, a player restricted to computable strategies, is allied to Nimue, the player choosing $P$, and both play against Merlin, the player choosing $v$. (The details of the game are recalled in \autoref{section-arthur-nimue-merlin-games}.)

While Arthur–Nimue–Merlin games were introduced recently by Kihara, they are a concrete description of objects that arose much earlier, in a different guise, in the study of realizability. To reassure the reader unfamiliar with realizability from the outset, this paper is meant to be accessible with only a background in basic computability theory, except for \autoref{section-topos-theory} where we explain how the proof of our main result can be understood from the point of view of realizability. Nevertheless, let us recall the historical link.

The effective topos, defined by Hyland, is a model of constructive mathematics — more precisely, an elementary topos with natural numbers object and thus a model of intuitionistic higher-order logic — which can be viewed as the “world of computable mathematics”. This topos contains the topos of sets, namely the “world of classical mathematics”, as a smallest non-degenerate subtopos, and in between them is an array of subtoposes, which are therefore intermediate worlds between computable and classical mathematics.

In his seminal work introducing the effective topos, Hyland associated to each Turing degree $\mathbf{d}$ a subtopos of “$\mathbf{d}$-computable mathematics”. This subtopos admits a description similar to the effective topos with plain computability replaced by computability relative to $\mathbf{d}$ (this was checked by Phoa \cite{Phoa}). Hyland proved that this defines an embedding of Turing degrees into subtoposes of the effective topos \cite[Theorem 17.2]{HylandEff}. However, there are many subtoposes which do not arise in this way.

As a general fact of topos theory, the subtoposes of an elementary topos are described by the so-called Lawvere–Tierney topologies (also known as modalities, local operators, or nuclei) on the topos. This characterization is already much more concrete: Lawvere–Tierney topologies on the effective topos can be defined in a simple way without reference to the effective topos, as maps $\P(\N) \to \P(\N)$ satisfying a few reasonable properties to expect of an oracle (see the beginning of \autoref{section-topos-theory}). Nevertheless, it remains somewhat mysterious. A first step towards a more explicit description was undertaken by Lee and van Oosten \cite{LeeVanOosten}, who introduced the subclass of \emph{basic} Lawvere–Tierney topologies on the effective topos, and described them using certain well-founded trees called \emph{sights}. However, the main advance was made by Kihara \cite{KiharaLT}, who proposed Arthur–Nimue–Merlin degrees as a considerably simpler description which encompasses all Lawvere–Tierney topologies on the effective topos. This description also sheds light on Hyland's embedding of Turing degrees into Arthur–Nimue–Merlin degrees: it is obtained by identifying a Turing oracle $g : \N \to \N$ with the Arthur–Nimue–Merlin oracle $\hat{g} : \N \to \P(\P(\N)), n \mapsto \{\{g(n)\}\}$ (the choices of Nimue and Merlin are forced, so that Arthur is in essence using an ordinary, deterministic Turing oracle). Lee and van Oosten's basic degrees, on the other hand, are those represented by oracles which are constant functions $\N \to \P(\P(\N))$.

This provides a first motivation for studying Arthur–Nimue–Merlin degrees, among the many kinds of degrees that appeared in computability theory in the last 70 years (such as enumeration degrees, Medvedev degrees, Muchnik degrees, etc.): although the details of the Arthur–Nimue–Merlin game may at first appear arbitrary, in fact it arises as a more explicit description of all “reasonable” oracles, which are also all “intermediate worlds between computable and classical mathematics”.

Another motivation is that Arthur–Nimue–Merlin degrees seem to be a natural completion of Turing degrees in the sense that several notions studied in classical computability theory find a natural home in the Arthur–Nimue–Merlin degrees. For example, define $g(e) := \{\{0\}\}$ if the $e$-th program with no input and boolean output terminates and returns false, $g(e) := \{\{1\}\}$ if it terminates and returns true, and $g(e) := \{\{0, 1\}\}$ if it does not terminate. This corresponds informally to having an oracle that will give Arthur the return value of a boolean-valued program that terminates, but may return either boolean for one that does not terminate.\footnote{This degree is called $\mathtt{LLPO}$ by Kihara \cite[Example 2.8 in published version, 2.5 in arXiv version]{KiharaLT} due to its connection with the lesser limited principle of omniscience in constructive mathematics.} The Turing degrees above $g$ are exactly the so-called “PA degrees”, a classical object of study in computability (see \cite[Definition 15.3.11]{CooperBook} and \cite[Theorem 10.3.3]{Soare}). Furthermore, $g$ is an intermediate degree between the computable degree $\mathbf{0}$ and the degree of the halting problem $\mathbf{0}'$, while, famously, no natural such degree is known in the Turing degrees.

Another example is the oracle $g$ where $g(n)$ does not depend on $n$ and consists of the cofinite subsets of $\N$. Informally, this means that Nimue chooses a natural number, then Merlin can replace it by any larger number before Arthur can see it. This example is due to Pitts \cite[Example 5.8]{PittsThesis}, and van Oosten proved \cite[Theorem 2.2]{vanOostenPittsOperator} that the functions $\N \to \N$ computable from this oracle exactly coincide with the hyperarithmetical functions. These can be defined in many other ways (e.g., as the $\Delta_1^1$ functions in the analytical hierarchy, or as the functions belonging to the $\omega_1^{\operatorname{CK}}$-th level of the constructible hierarchy), and have been intensely studied in effective descriptive set theory, where they are the analogue of Borel functions in classical descriptive set theory.

Recently, Kihara and (Ming) Ng studied the preorder on basic Arthur–Nimue–Merlin degrees defined by $\mathbf{a} \leq^{\mathsf{o}} \mathbf{b}$ when $\mathbf{a} \leq \mathbf{b} \vee \mathbf{t}$ for some Turing degree $\mathbf{t}$, exhibiting a connection with the Katětov order in set theory \cite{KiharaNg2,KiharaNg2Bis}, while Ahman and Bauer explained how to understand the characterization of Lawvere–Tierney topologies on the effective topos in terms of Arthur–Nimue–Merlin games in a more conceptual way via “oracle modalities” \cite{AhmanBauer}, which we use in \autoref{section-topos-theory}.

\subparagraph*{Contribution}
We make a further advance in the understanding of Arthur–Nimue–Merlin degrees by exhibiting an embedding of the Turing degrees with the opposite of their standard order into the Arthur–Nimue–Merlin degrees — in fact into the basic ones. Explicitly, this means that we define for each Turing degree $\mathbf{a}$ a basic Arthur–Nimue–Merlin degree $\coT(\mathbf{a})$, called its \emph{co-Turing degree}, so that the map $\coT$ is an order-reversing embedding, i.e., $\mathbf{a} \leq \mathbf{b} \; \iff \; \coT(\mathbf{a}) \geq \coT(\mathbf{b})$ for all Turing degrees $\mathbf{a}, \mathbf{b}$. Furthermore, we show that a non-zero Turing degree is never comparable with the co-Turing degree of a non-zero Turing degree.

Our “co-Turing degree” construction provides one answer to \cite[section 6, Question 1]{KiharaLT} which asked for other “interesting” basic Arthur–Nimue–Merlin degrees.

\subparagraph*{Outline}
\autoref{section-arthur-nimue-merlin-games}
precisely describes the Arthur–Nimue–Merlin game and recalls their basic properties. This is essentially a summary of some parts of \cite{KiharaLT} but with a few inessential differences in notation and terminology. In \autoref{section-decoding-and-coding-degrees} we define the \emph{decoding degree} and \emph{coding degree} associated with a promise problem and we prove our main technical result, \autoref{main-strategy-result}. In \autoref{section-co-turing-degrees}, we define co-Turing degrees as a particular case of coding degrees and prove our main theorem, \autoref{reverse-embedding-theorem}, according to which the mapping from Turing degrees to co-Turing degrees is an order-reversing embedding, as well as \autoref{incomparability-turing-and-co-turing-degrees} on the incomparability of Turing and co-Turing degrees. \autoref{section-topos-theory} is devoted to reinterpreting the construction of co-Turing degrees in terms of Lawvere–Tierney topologies on the effective topos and is the only section which requires familiarity with the effective topos (after an initial recapitulation on Lawvere–Tierney topologies on the effective topos, which is self-contained).

\subparagraph*{Notations}
We work in classical logic unless otherwise specified. However, we still adopt the convention from constructive mathematics to call \textbf{inhabited} a set which has an element (this is classically equivalent to being non-empty). We write $f \colon A \to B$ for a total function and $f \colon A \dasharrow B$ for a partial function. In the latter case, we write $f(x){\downarrow}$ to mean that $f$ is defined at $x$, and $f(x){\downarrow} = y$ to mean that furthermore the value equals $y$. The powerset of $X$ is denoted by $\P(X)$. We fix a computable Gödel coding of tuples by natural numbers, denoted $(m_0, \dots, m_{k-1}) \mapsto \langle m_0, \dots, m_{k-1} \rangle$. The $e$-th partial computable function is $\varphi_e \colon \N \dasharrow \N$ and the $e$-th partial computable function with oracle $f \colon \N \to \N$ is $\varphi_e^f$. (The numbering of programs is independent of the oracles, e.g., there exists $e$ such that $\varphi_e^f = f$ for all $f$.)

\section{Arthur–Nimue–Merlin games and degrees}\label{section-arthur-nimue-merlin-games}

A \textbf{Turing oracle} is a function $\N \to \N$. Turing oracles are preordered by Turing reducibility $\preceq$, and the partially ordered set of Turing degrees $\DT$ is defined as the quotient of Turing oracles by the canonical associated equivalence relation, $f \equiv g$ when $f \preceq g$ and $g \preceq f$, equipped with the partial order $\leq$ induced by $\preceq$. We denote the degree of an oracle $f$ by $[f]$.

The story is entirely similar for Arthur–Nimue–Merlin degrees, but the reducibility is more complicated.

\begin{definition}\label{definition-ANM}
An \textbf{Arthur–Nimue–Merlin oracle} is a function $\N \to \P(\P(\N))$.

Let $f, g$ be Arthur–Nimue–Merlin oracles. The \textbf{reduction game} from $f$ to $g$ is a game between three players named Arthur, Nimue and Merlin, played as follows. Arthur and Nimue are allied against Merlin. The game starts with an initial configuration, which consists of a natural number $n$ and a target set $T \in f(n)$. Nimue and Merlin can see both $n$ and $T$, while Arthur can only see $n$. Accordingly, we call $n$ the \emph{public part} of the initial configuration, and $T$ the \emph{secret part}. The moves are played in the following order: Arthur, then Nimue, then Merlin play in turn, repeatedly, until Arthur, during his turn, plays a special move that ends the game (or until a player cannot play and loses). When it is his turn to play, Arthur has two kinds of moves available: either declare a natural number of his choice as the “final value”, which ends the game, or query the oracle at a natural number of his choice, the “question”. If Arthur declares $t \in \N$ as the final value, then Arthur and Nimue win if $t \in T$ (where $T \subseteq \N$ is the secret part of the initial configuration), otherwise Merlin wins. If Arthur queries the oracle on a question $q \in \N$, both Nimue and Merlin can see the question, and it is Nimue's turn to play next. Nimue chooses a set $R \in g(q)$. (If $g(q)$ is empty, then she has no move available, so she and Arthur lose immediately.) Nimue communicates the choice of $R$ to Merlin, but Arthur cannot see it. Then Merlin chooses a reply $r \in R$ (if $R$ is empty, Merlin has no move available and immediately loses). Both Arthur and Nimue can see $r$, and the game continues with Arthur's turn. If the game never ends, then Merlin wins by default.

A \textbf{Nimue strategy} for this game is a function taking an initial configuration and a finite sequence of moves by Merlin and Arthur to an element of $\mathcal{P}(\mathbb{N})$ to be played by Nimue. An \textbf{Arthur strategy} is a \emph{partial computable} function taking a natural number (representing the public part, $n$, of the initial configuration and a finite sequence of natural numbers (representing moves by Merlin) to a move to be played by Arthur (coded, e.g., as either $\langle 0, q \rangle$ to query the oracle for $n$, or $\langle 1, t \rangle$ to declare $t$ as the final value). (If the move returned by the strategy is undefined or not of the appropriate form, then Arthur and Nimue will lose the game.)

An \textbf{Arthur–Nimue strategy} is a pair consisting of an Arthur strategy and a Nimue strategy. It is \textbf{winning} when Arthur and Nimue win the game when they play following that strategy, no matter what the initial configuration and Merlin's moves are.\footnote{We do not define the notion of winning Merlin strategy because Merlin cannot have one except in trivial cases: Arthur might win “by accident” by declaring a good value on his very first move.} The oracle $f$ is \textbf{reducible} to $g$, written $f \preceq g$, when such a winning strategy exists.
\end{definition}

This definition is due to Kihara \cite[Definitions 2.13 \& 2.14\footnote{Numbered 2.9 and 2.10 in the arXiv version}]{KiharaLT}. Our description differs from Kihara's in two inessential ways. The first is that for technical convenience, we used functions $g : \N \to \P(\P(\N))$ for the oracles instead of what is called a “bilayer function” in \cite{KiharaLT}, namely a partial function $g : \N \times \Lambda \dasharrow \P(\N)$, where $\Lambda$ is an unspecified set that varies from oracle to oracle (the idea being that Arthur chooses a question $q$, Nimue chooses a “parameter” $\lambda$, they must ensure that $g(n, \lambda)$ is defined, and Merlin returns an element of this set). The second difference is that we consider the data $(n, T)$ as an “initial configuration” rather than an opening move by Merlin, which is immaterial but more convenient when we need to speak of the sequence of moves by Merlin.

\begin{proposition}[Kihara {\cite[Proposition 2.19\footnote{Numbered 2.13 in the arXiv version}]{KiharaLT}}]
The reducibility relation $\preceq$ is a preorder on Arthur–Nimue–Merlin oracles.
\end{proposition}
\begin{proof}
We leave reflexivity as an exercise. Assume $f \preceq g$ through an Arthur strategy $s$ and a Nimue strategy $S$, and $g \preceq h$ through $s'$ and $S'$. Then $f \preceq h$ through the following strategy. Arthur receives $n$ and must find an element of some $T \in f(n)$. He simulates $s$ until it makes a query $q$. Nimue simulates $S$ to choose $R \in g(q)$. Then Arthur simulates $s'$ using $q$ as the public part of the initial configuration, and Nimue simulates $S'$ using $q$ and $R$ as the initial configuration. As $s'$ and $S'$ reduce $g$ to $h$, this ultimately allows Arthur to compute an element $r \in R$. Then Arthur and Nimue pretend $r$ was played by Merlin and go back to simulating $s$ and $S$. They continue in the same way, simulating each oracle query via a full run of $s'$ and $S'$, until $s$ declares a final value $t \in T$ which Arthur returns.
\end{proof}

\begin{definition}\label{anm-degrees}
Two Arthur–Nimue–Merlin oracles $f, g : \N \to \P(\P(\N))$ are \textbf{equivalent}, denoted $f \equiv g$, when $f \preceq g$ and $g \preceq f$. The partially ordered set of \textbf{Arthur–Nimue–Merlin degrees}, denoted $\DANM$ is the quotient of Arthur–Nimue–Merlin oracles by equivalence, equipped with the order $\leq$ induced by $\preceq$. We write $[f]$ for the degree of an oracle $f$.
\end{definition}

As announced, Arthur–Nimue–Merlin degrees generalize Turing degrees.

\begin{remark}
Let $f, g : \N \to \N$ be Turing oracles. Define the Arthur–Nimue–Merlin oracles $\hat{f} : \N \to \P(\P(\N))$ by $\hat{f}(n) := \{\{f(n)\}\}$ and likewise $\hat{g}$. Then $f \preceq g$ (i.e., $f$ is Turing-reducible to $g$) if and only if $\hat{f} \preceq \hat{g}$ (i.e., $\hat{f}$ is Arthur–Nimue–Merlin reducible to $\hat{g}$). In particular, $[f] \mapsto [\hat{f}]$ is a well-defined embedding of the partially ordered set of Turing degrees $\DT$ into the partially ordered set of Arthur–Nimue–Merlin degrees $\DANM$.
\end{remark}

In the sequel, we implicitly identify a Turing oracle $f$ with the Arthur–Nimue–Merlin oracle $\hat{f}$ and we view Turing degrees as a subset of Arthur–Nimue–Merlin degrees.

The \textbf{bottom degree} $\mathbf{0}$ is the Turing degree of computable functions. It is the least Arthur–Nimue–Merlin degree since Arthur can simply compute the function without querying the oracle at all.

We also single out the basic degrees in the sense of Lee and van Oosten \cite{LeeVanOosten}.

\begin{definition}
A \textbf{basic oracle} is a subset $G \subseteq \P(\N)$, which we implicitly identify with the Arthur–Nimue–Merlin oracle defined by the constant function $n \mapsto G$. A \textbf{basic degree} is a degree represented by some basic oracle.
\end{definition}

In the reduction game to a basic oracle, the questions asked by Arthur are irrelevant. Evocatively, Nimue tries to communicate information to Arthur, but with Merlin scrambling it on the way, and Arthur can listen for as long as he wants but must ultimately decide on an answer to return.

Conversely, in the reduction game from a basic oracle, Arthur starts with no information at all on the target set of which he must declare an element, and must entirely rely on the help of Nimue to obtain such information.

\begin{remark}\label{remark-when-basic-oracle-is-zero}
In particular, a basic oracle $G \subseteq \P(\N)$ has degree $\mathbf{0}$ if and only if $\bigcap G$ is inhabited (this is \cite[Proposition 3.1(iii)]{LeeVanOosten}) since, in the reduction game from $G$ to a computable oracle, there is no help available, so Arthur can only win by always returning the same constant.
\end{remark}

There is no greatest Turing degree since every Turing degree $\mathbf{d}$ is less than its Turing jump $\mathbf{d}'$. On the other hand, there is a “degenerate” greatest Arthur–Nimue–Merlin degree and a greatest element of the rest.

\begin{definition}\label{bottom-top-and-omniscient-degrees}
The \textbf{top} degree $\top \in \DANM$ is the basic degree represented by the basic oracle $\{\varnothing\}$. The \textbf{omniscient} degree\footnote{The omniscient degree is denoted $\mathtt{Error}_{1/2}$ by Kihara in \cite{KiharaLT}.} $\Om \in \DANM$ is the basic degree represented by the basic oracle $\{\{k\} \mid k \in \N\}$.
\end{definition}

\begin{remark}
The top degree $\top$ is the greatest Arthur–Nimue–Merlin degree and the omniscient degree $\Om$ is the greatest Arthur–Nimue–Merlin degree different from the top degree.
\end{remark}

\begin{proof}
If $f$ is an oracle such that $\varnothing \in f(n)$ for some $n$, then Arthur and Nimue can win the reduction game from any oracle to $f$ through the following strategy: Arthur queries the oracle on $n$ and Nimue plays $\varnothing$, so that Merlin loses by lack of available move. In particular, $\top$ is the maximum Arthur–Nimue–Merlin degree.

Now let $f$ be an oracle whose degree is not $\top$. Then $f$ reduces to $\Om$ by the following Arthur–Nimue strategy. Consider an initial configuration consisting of $n \in \N$ and $T \in f(n)$. By the previous remark, $T$ is inhabited. Arthur and Nimue can win the reduction game from $f$ to $\Om$ by the following strategy: Arthur queries the oracle on $n$; Nimue picks any $t \in T$ and plays $\{t\}$; Merlin is forced to answer $t$ to Arthur, and Arthur returns $t$.
\end{proof}

\begin{remark}
The basic oracle $\{\{0\}, \{1\}\}$ has degree $\Om$ (if Nimue can repeatedly communicate one bit to Arthur, she can communicate an arbitrary natural number, by encoding it in unary).
\end{remark}

Finally, we note that least upper bounds exist in the Arthur–Nimue–Merlin degrees and are just as easy to compute as in the Turing degrees.

\begin{remark}
Let $f$ and $g$ be Arthur–Nimue–Merlin oracles. The least upper bound of their degrees exists and can be explicitly constructed as $[f] \vee [g] = [h]$ where $h(2n) := f(n)$ and $h(2n+1) := g(n)$. The inclusion of Turing degrees into Arthur–Nimue–Merlin degrees preserves least upper bounds.
\end{remark}

\section{Decoding and coding degrees for promise problems}\label{section-decoding-and-coding-degrees}

We borrow terminology from complexity theory \cite{GoldreichPromiseProblems} by calling \textbf{promise problem} a pair $(P, Q)$ of disjoint inhabited subsets $P, Q \subseteq \N$, informally understood as the computational problem “given an element $n \in P \cup Q$, determine whether $n \in P$ or $n \in Q$”.\footnote{This should not be confused with the \emph{separation problem} $(P, Q)$, which requires the program solving the problem to terminate (with any value) even on inputs outside $P \cup Q$.} (Requiring $P$ and $Q$ to be inhabited excludes annoying special cases in what follows.)

\begin{definition}\label{definition-coding-decoding-degrees}
Let $(P, Q)$ be a promise problem. The \textbf{decoding oracle} $d_{P, Q}$ is defined by $d_{P, Q}(n) = \{\{0\}\}$ for $n \in P$ and $d_{P, Q}(n) = \{\{1\}\}$ for $n \in Q$ and $d_{P, Q}(n) = \varnothing$ for $n \notin P \cup Q$. The \textbf{coding oracle} $c_{P, Q}$ is the basic oracle defined by $c_{P, Q}(n) := \{P, Q\}$. At the level of degrees, we write $\mathbf{d}_{P, Q} := [d_{P, Q}]$ and $\mathbf{c}_{P, Q} := [c_{P, Q}]$.
\end{definition}

The decoding oracle $d_{P, Q}$ is the most natural Arthur–Nimue–Merlin oracle associated to the promise problem $(P, Q)$. It is decidable (has degree $\mathbf{0}$) precisely when $(P, Q)$ is decidable as a promise problem in the sense just described: given $n \in P \cup Q$, it can be computably determined whether $n \in P$ or $n \in Q$. As an oracle, it answers questions of the form “is $n$ in $P$ or in $Q$?”, where questions $n \notin P \cup Q$ are forbidden. Remark that when $P$ and $Q$ are complements of each other, $\mathbf{d}_{P, Q}$ is the ordinary Turing degree of both.

The coding oracle is more interesting: it communicates bits helpfully chosen by Nimue in a form adversarially encoded by Merlin as instances of the promise problem $(P, Q)$. Namely, Merlin returns an element of $P$ if Nimue chose $0$ and an element of $Q$ if Nimue chose $1$.

\begin{remark}\label{decoding-plus-coding-is-omniscient}
We have $\mathbf{d}_{P, Q} \vee \mathbf{c}_{P, Q} = \Om$ because this oracle allows both to obtain answers to arbitrary questions coded as instances of the promise problem $(P, Q)$ and to decode these answers (and it is not $\top$ because we required $P$ and $Q$ to be inhabited).
\end{remark}

It is intuitively plausible that the harder the promise problem $(P, Q)$ is, the weaker the $c_{P, Q}$ oracle is because it becomes harder for Arthur to decode the messages sent by Nimue. Thus, the higher $\mathbf{d}_{P, Q}$, the lower $\mathbf{c}_{P, Q}$. However, this is not strictly true.

\begin{example}
In our opening riddle, we asked whether $\mathbf{c}_{K, K^c} = \Om$ where $K := \{e \mid \varphi_e(0){\downarrow}\}$ is the halting set and $K^c$ its complement. We observed that this is indeed the case, because Nimue can transmit a bit to Arthur by either playing $K$ followed by $K^c$, or $K^c$ followed by $K$. Arthur receives two Turing machines and runs them in parallel, so that which one terminates reveals the bit chosen by Nimue. This shows that $\mathbf{d}_{P, Q}$ can be made non-zero while $\mathbf{c}_{P, Q}$ remains omniscient.
\end{example}

More generally, if $\mathbf{d}_{P \times Q, \, Q \times P} = \mathbf{0}$,\footnote{The notation suppresses the implicit conversion of pairs of natural numbers to natural numbers through the pairing $\langle -, - \rangle$.} then $\mathbf{c}_{P, Q} = \Om$ because the following Arthur–Nimue strategy allows Nimue to communicate a bit to Arthur: he queries the oracle twice, and depending on the bit, she plays $P$ then $Q$ or $Q$ then $P$, and Arthur solves the resulting instance of the promise problem $(P \times Q, \, Q \times P)$.

\begin{example}\label{example-dedekind-cuts}
In fact, $\mathbf{d}_{P, Q}$ can be made an arbitrary Turing degree while keeping $\mathbf{c}_{P, Q}$ omniscient. Let $(P, Q)$ be the Dedekind cut of a real number $x$, namely $P := \{q \in \Q \mid q < x\}$ and $Q := \{q \in \Q \mid q > x\}$ (viewed as subsets of $\N$ through a Gödel coding of rational numbers). Then $\mathbf{d}_{P, Q}$ is the Turing degree of $x$ in the standard sense, but $(P \times Q, \, Q \times P)$ can be decided by simply comparing rationals.
\end{example}

\begin{example}
For another example of the same phenomenon, let $f : \N \to \{0, 1\}$ and define $P$ as the set of codings of a prefix of this function, $\langle f(0), \dots, f(n-1) \rangle$, such that $f(n) = 0$, and $Q$ likewise with $f(n) = 1$. Then $\mathbf{d}_{P, Q} = [f]$ (for the $\geq$ direction, compute $f(0), f(1), \dots$ iteratively using the oracle until the desired value), yet $(P \times Q, \, Q \times P)$ is decidable because if we have a prefix of $f$ in $P$ and another prefix in $Q$, one of them is strictly longer and reveals which of $P$ and $Q$ the other belongs to.
\end{example}

We now arrive at our main technical result: the implication $\mathbf{d}_{P \times Q, \, Q \times P} = \mathbf{0} \; \implies \; \mathbf{c}_{P, Q} = \Om$ is in fact an equivalence. Moreover, the result relativizes to an arbitrary Turing degree, and this fact will be the crux of the proof that the construction of co-Turing degrees defines an embedding.

\begin{theorem}\label{main-strategy-result}
\textit{(Plain version.)} Let $(P, Q)$ be a promise problem. We have $\mathbf{c}_{P, Q} = \Om$ if and only if $\mathbf{d}_{P \times Q, \, Q \times P} = \mathbf{0}$.

\textit{(Relativized version.)} Let $(P, Q)$ be a promise problem and let $\mathbf{a}$ be a Turing degree. We have $\mathbf{c}_{P, Q} \vee \mathbf{a} = \Om$ if and only if $\mathbf{d}_{P \times Q, \, Q \times P} \leq \mathbf{a}$.
\end{theorem}

\begin{proof}
Let us first prove the “plain” version. For the entirety of the proof, we fix elements $p \in P$ and $q \in Q$, which exist by assumption.

We assume that $\mathbf{c}_{P, Q}$ is omniscient, which means that
there exists a winning Arthur–Nimue strategy in the game reducing $\Om$
to $\mathbf{c}_{P, Q}$.
We fix such a strategy and call $\sigma$ Arthur's part of the strategy;
thus, $\sigma$ takes a finite sequence of natural numbers representing
moves by Merlin to a move by Arthur, which can be “query the oracle
once more” or “declare a final value $u$”.

Given an element $(x,y) \in (P\times Q) \,\cup\, (Q\times P)$, we consider
the function ${M_{(x,y)} \colon \{\Psymb,\Qsymb\}^2 \to
\mathbb{N}}$ (where $\Psymb,\Qsymb$ are two
labels\footnote{If the reader has difficulty distinguishing
 $\Psymb,\Qsymb$ from $P,Q$ typographically, we hasten
 to reassure that this will not cause any confusion.})
defined explicitly by
\begin{align*}
M_{(x,y)}(\Psymb,\Psymb) &= p,
\\ M_{(x,y)}(\Psymb,\Qsymb) &= x,
\\ M_{(x,y)}(\Qsymb,\Psymb) &= y,
\\ M_{(x,y)}(\Qsymb,\Qsymb) &= q.
\end{align*}
This function is designed to ensure the following property (\textdagger):
if $(x,y) \in P\times Q$, then $M_{(x,y)}(\Psymb, r) \in P$ and
$M_{(x,y)}(\Qsymb, r) \in Q$ for any $r \in
\{\Psymb,\Qsymb\}$, and symmetrically, if $(x,y) \in Q\times P$, then
$M_{(x,y)}(\ell, \Psymb) \in P$ and $M_{(x,y)}(\ell, \Qsymb)
\in Q$ for any $\ell \in \{\Psymb,\Qsymb\}$.

We define an ancillary game, which is a two-player imperfect
information game parameterized by $(x, y) \in (P \times Q) \,\cup\, (Q \times P)$.
The two players are “Left” and “Right”. On the $i$-th turn,
Left chooses a value $\ell_i \in \{\Psymb, \Qsymb\}$, and
Right simultaneously chooses $r_i \in \{\Psymb, \Qsymb\}$. (Each player
has seen the other's moves up to this turn, but the two moves on this
turn are chosen independently.) What happens next depends on the
value of Arthur's strategy $\sigma$ on the sequence
$M_{(x, y)}(\ell_0, r_0), \ldots, M_{(x, y)}(\ell_i, r_i)$.
If $\sigma$ returns “query the oracle”, then the ancillary game continues.
If $\sigma$ returns “declare a final value $u$” where $u \in \{0, 1\}$,
then Left wins if $u = 0$, whereas Right wins if $u = 1$. If $\sigma$
fails to return or returns a malformed value, or if the game lasts
infinitely long, then the result is a draw (no player wins).

Assume temporarily that $(x, y) \in P \times Q$.
We observe that Left has a winning strategy in the ancillary game.
Indeed, by (\textdagger), Left can control whether the successive
elements provided to $\sigma$ are in $P$ or $Q$, and using Nimue's
part of the Arthur–Nimue winning strategy, this enables Left to
make Arthur's strategy output 0.
Let $L$ be a winning strategy for Left. This is a function from finite
sequences in $\{\Psymb, \Qsymb\}$, representing past moves by Right,
to elements of $\{\Psymb, \Qsymb\}$, representing the next move by Left.
For every sequence $s \in \{\Psymb, \Qsymb\}^\N$ of moves by Right, the
play where Left plays according to $L$ and Right plays the successive
elements of $s$ must end in a victory for Left after a finite number $m(s)$
of moves. The function $m : \{\Psymb, \Qsymb\}^\N \to \N$ is
continuous (where $\{\Psymb, \Qsymb\}^\N$ is topologized by the product
of the discrete topology)
because $m(s)$ only depends on the first $m(s)$ terms of $s$.
By compactness of the Cantor space $\{\Psymb, \Qsymb\}^\N$, it
must be bounded by some constant $M$. This means that a finite fragment
of $L$ — its restriction to finite sequences in $\{\Psymb, \Qsymb\}$
of length at most $M$ — already suffices to witness that Left has
a winning strategy.

We have proved that for all $(x, y) \in P \times Q$, there exists
a finite fragment of strategy, namely a function only defined on
finite sequences in $\{\Psymb, \Qsymb\}$ up to a certain length $M$,
which certifies that Left can win, in the sense that running the
game up to the $M$-th turn using this strategy fragment for the moves
of Left and any sequence of length $M$ in $\{\Psymb, \Qsymb\}$
for the moves of Right results in the game already ending with a
victory for Left.
 Symmetrically, for $(x, y) \in Q \times P$, there exists a finite
strategy fragment which certifies that Right can win. Finally,
the two conditions are mutually exclusive since the players cannot
both win at the same time. This shows that $(x, y) \in P \times Q$
if and only if Left has a winning finite strategy fragment and
$(x, y) \in Q \times P$ if and only if Right has a winning finite strategy
fragment.

We deduce an algorithm to decide, given $(x, y) \in (P \times Q)
\,\cup\, (Q \times P)$, whether $(x, y) \in P \times Q$
or $(x, y) \in Q \times P$: enumerate finite strategy fragments in parallel
to search at the same time for a winning strategy fragment for Left
and a winning strategy fragment for Right. Whether a strategy
fragment is winning for a given player can be semi-decided by
simulating the plays with this fragment for the given player and
all possible sequences up to the length of the fragment for the other
player.

This concludes the proof of the “plain” version of the theorem.
To prove the “relativized” version, we make the following observation:
the proof of the “plain” version which we have just given applies just
as well if we replace computability by computability relative to a
Turing oracle $A \subseteq \mathbb{N}$. The assumption that $A$ is
a Turing oracle and not a general Arthur–Nimue–Merlin oracle is used
so that the ancillary game is well-defined as Arthur's strategy $\sigma$
is using a deterministic oracle, as well as to ensure it is possible
to run several computations in parallel under this oracle, some of which
are not guaranteed to be well-behaved. (This might be impossible under
an oracle $g$ such that some of the sets $g(q)$ are empty, which means
that $q$ is a “forbidden question” that makes the computation immediately
fail if asked to the oracle, including any other “threads” that the program
might be running.)
This gives us the following result: if Arthur–Nimue have a winning strategy in the
game reducing $\Om$ to $\mathbf{c}_{P,Q}$, but where now Arthur is
allowed to access the oracle $A$ in computing his strategy, then there
is an $A$-computable way to solve to the promise problem $(P\times Q,
\, Q\times P)$. The hypothesis is equivalent to Arthur–Nimue having
a winning strategy in the game reducing $\Om$ to
$[A] \vee \mathbf{c}_{P,Q}$ (because an
$A$-computable strategy for Arthur is the same thing as a strategy for
Arthur that can make the move of querying the $A$-oracle), and
the conclusion is equivalent to Arthur having a winning strategy
in the game reducing $d_{P\times Q, \, Q\times P}$ to $A$.
\end{proof}

\begin{remark}\label{remark-answer-to-riddle}
This provides an answer to the second part of our opening riddle, which asks whether $\mathbf{c}_{A, A^c} = \Om$ where $A := (K \times K) \cup (K^c \times K^c)$ (here $K$ is the halting set), and thus $A^c = (K \times K^c) \cup (K^c \times K)$. If such were the case, then by \autoref{main-strategy-result} the promise problem $(A \times A^c, \, A^c \times A)$ would be decidable. This would allow to decide whether a program $e$ halts, by deciding the instance $((\mathtt{loop}, e), (\mathtt{halt}, e))$ of this problem, where $\mathtt{loop}$ is a fixed program that does not halt and $\mathtt{halt}$ is a fixed program that does.
\end{remark}

\section{Co-Turing degrees}\label{section-co-turing-degrees}

We turn to the construction of co-Turing degrees, as special cases of $\mathbf{c}_{P, Q}$. For technical convenience, we use Turing degrees of subsets $A \subseteq \N$ (namely the Turing degrees of their characteristic functions $\N \to \{0, 1\}$) instead of functions $\N \to \N$.

One might naïvely attempt to define the co-Turing degree of $A \subseteq \N$ to be $\mathbf{c}_{A, A^c}$, but this is not a function of the Turing degree of $A$. For example, we have just seen (\autoref{remark-answer-to-riddle}) that there is an undecidable set $A$ for which $\mathbf{c}_{A, A^c} < \Om$, and using Dedekind cuts (\autoref{example-dedekind-cuts}), we can find another set $B$ with the same Turing degree such that $\mathbf{c}_{B, B^c} = \Om$. Instead, we use a slight modification, for which a conceptual explanation will be provided in \autoref{section-topos-theory}.

\begin{definition}\label{definition-codegrees}
Let $A \subseteq \N$. For $i \in \{0, 1\}$, we let $H_i^A$ be the set of programs without argument which return $i$ when run with $A$ as an oracle (i.e., $H_i^A := \{e \mid \varphi_e^A(0){\downarrow} = i\}$). We define the \textbf{co-Turing oracle} of $A$ as the basic oracle $\coT(A) := c_{H_0^A, H_1^A}$.
\end{definition}

In words, the oracle $\coT(A)$ cooperatively answers arbitrary boolean questions but, instead of directly providing the answers, it adversarially encodes each answer in the form of a program which returns the answer when run with $A$ as an oracle.

\begin{theorem}
The map $\coT$ from Turing oracles with the opposite of their normal preorder to Arthur–Nimue–Merlin oracles (with their normal preorder) is order-preserving and order-reflecting.
\end{theorem}

\begin{proof}
Let $A, B \subseteq \N$.

Assume $A \preceq B$. We seek to prove $\coT(B) \preceq \coT(A)$. Unfolding the reduction game, Arthur needs to declare a program with oracle $B$ which outputs a bit known to Nimue but not to him. The following Arthur–Nimue strategy works: Arthur queries the oracle once; Nimue chooses this bit, so that Merlin gives Arthur a program with oracle $A$ which outputs this bit; then Arthur modifies this program to replace each call to the oracle $A$ by a full run of a fixed Turing reduction from $A$ to $B$.

Conversely, assume $\coT(B) \preceq \coT(A)$ and let us prove $A \preceq B$. We have $\coT(B) \vee [B] = \Om$ because $\coT(B) \vee \mathbf{d}_{H_0^B, H_1^B} = \Om$ (\autoref{decoding-plus-coding-is-omniscient}) and $\mathbf{d}_{H_0^B, H_1^B} \leq [B]$. Since the join is order-preserving in both arguments, $\coT(A) \vee [B] = \Om$. By \autoref{main-strategy-result}, this entails $d_{H_0^A \times H_1^A, H_1^A \times H_0^A} \preceq B$. It remains to see that $A \preceq d_{H_0^A \times H_1^A, H_1^A \times H_0^A}$ (this is in fact an equality since the converse inequality is trivial). We use the following reduction. Given $n \in \N$, we construct a program $r_0$ with oracle $A$ which returns $0$ if $n \in A$ and $1$ if $n \notin A$, and likewise $r_1$ which returns $1$ if $n \in A$ and $0$ if $n \notin A$. Then $(r_0, r_1) \in H_0^A \times H_1^A$ if $n \in A$ and $(r_0, r_1) \in H_1^A \times H_0^A$ if $n \notin A$ so querying the $d_{H_0^A \times H_1^A, H_1^A \times H_0^A}$ oracle on $(r_0, r_1)$ solves the problem.
\end{proof}

\begin{corollary}\label{reverse-embedding-theorem}
The map $\coT$ induces an embedding $\coT : (\DT)^{\op} \inject \DANM$ of the partially ordered set of Turing degrees with the opposite of their normal order into the partially ordered set of Arthur–Nimue–Merlin degrees.
\end{corollary}

We call $\coT(\mathbf{a})$ the \textbf{co-Turing degree} of $\mathbf{a}$. The definition means that $\coT(\mathbf{a}) = [\coT(A)]$ for any subset $A \subseteq \N$ of Turing degree $\mathbf{a}$.

Having completely described the ordering inside the co-Turing degrees, we turn to the ordering between Turing and co-Turing degrees inside the Arthur–Nimue–Merlin degrees. The bottom degree $\mathbf{0}$ is a Turing degree below all co-Turing degrees, and the omniscient degree $\Om$ is a co-Turing degree (namely $\Om = \coT(\mathbf{0})$) above all Turing degrees. It turns out that apart from these two trivial cases, a Turing degree and a co-Turing degree are \emph{never} comparable. The proof involves the following lemma, which is a reformulation of a result of Kihara and (Keng) Ng \cite[Corollary 4.8]{KiharaNg1}, but we give a more direct proof, for the idea of which we are indebted to the pseudonymous MathOverflow user “Pastebee” \cite{MO507736}.

\begin{lemma}\label{lemma-T1-meet}
Let $A, B \subseteq \N$. If the promise problem $(H_0^A \times H_0^B, \, H_1^A \times H_1^B)$ is decidable, then one of $A$ and $B$ is decidable.
\end{lemma}

\begin{proof}
Let $e$ decide the promise problem $(H_0^A \times H_0^B, \, H_1^A \times H_1^B)$.

For $n \in \N$, let $r_0^n$ be the code of the oracle program that (ignores its argument and) returns $0$ if $n$ is in the oracle, otherwise $1$, and likewise let $r_1^n$ return $1$ if $n$ is in the oracle, otherwise $0$. We can arrange for $n \mapsto r_0^n$ and $n \mapsto r_1^n$ to be computable.

Consider the following putative algorithm to decide $A$. Given an input $m$, we use dovetailing to search in parallel for $n$ such that $\varphi_e(\langle r_1^m, r_1^n \rangle)$ and $\varphi_e(\langle r_1^m, r_0^n \rangle)$ are both defined and equal, and then return their common value.

Let us check the partial correctness of this algorithm. For a fixed $m$, if it finds such an $n$ in its search, then either the booleans $m \in A$ and $n \in B$ are the same, in which case $\varphi_e(\langle r_1^m, r_1^n \rangle)$ is this boolean, or $m \in A$ and $n \notin B$ are the same, in which case $\varphi_e(\langle r_1^m, r_0^n \rangle)$ is this boolean. Either way, the algorithm returns the boolean $m \in A$. Hence, if the algorithm always terminates, it decides $A$.

Now assume that $m_0$ is an input for which the above algorithm does not terminate. In this case, we claim that the following other algorithm decides $B$. Given an input $n$, we run two computations in parallel: in the first, if $\varphi_e(\langle r_1^{m_0}, r_1^n \rangle)$ terminates, we return it; in the second, if $\varphi_e(\langle r_1^{m_0}, r_0^n \rangle)$ terminates with a boolean output, we return its opposite boolean. The assumption that the previous algorithm does not terminate on $m_0$ means that $\varphi_e(\langle r_1^{m_0}, r_1^n \rangle)$ and the opposite of $\varphi_e(\langle r_1^{m_0}, r_0^n \rangle)$ can never be both defined and different, so it suffices to see that one of them must be defined and correct. If $m_0 \in A \iff n \in B$, the first is, while if $m_0 \in A \iff n \notin B$, the second is.
\end{proof}

\begin{theorem}\label{incomparability-turing-and-co-turing-degrees}
Let $\mathbf{a}$ and $\mathbf{b}$ be Turing degrees.
\begin{enumerate}
\item We always have $\coT(\mathbf{a}) \not\leq \mathbf{b}$.
\item We have $\coT(\mathbf{a}) \geq \mathbf{b}$ if and only if $\mathbf{a} = \mathbf{0}$ or $\mathbf{b} = \mathbf{0}$.
\end{enumerate}
\end{theorem}

\begin{proof}
\begin{enumerate}
\item Assume $\coT(\mathbf{a}) \leq \mathbf{b}$. From this and $\coT(\mathbf{a}) \vee \mathbf{a} = \Om$ (\autoref{decoding-plus-coding-is-omniscient}), we deduce ${\mathbf{b} \vee \mathbf{a} = \Om}$. We know that least upper bounds of Turing degrees in the Arthur–Nimue–Merlin degrees are computed as in the Turing degrees, so $\Om$ is a Turing degree, but it is not (no Turing oracle will let Arthur guess a secret number).
\item The “if” direction is clear. Assume $\coT(\mathbf{a}) \geq \mathbf{b}$. From this and $\coT(\mathbf{b}) \vee \mathbf{b} = \Om$ (\autoref{decoding-plus-coding-is-omniscient} again), we deduce $\coT(\mathbf{a}) \vee \coT(\mathbf{b}) = \Om$.

To put this in a form where \autoref{main-strategy-result} applies, we make the general observation that $\mathbf{c}_{P, Q} \vee \mathbf{c}_{P', Q'} = \mathbf{c}_{P \times P', Q \times Q'}$ for all $P, Q, P', Q' \subseteq \N$. (For the inequality $\leq$, Arthur is either asked to return a boolean known to Nimue as an instance of the promise problem $(P, Q)$, or as an instance of $(P', Q')$. In either case, he can obtain it as an instance of $(P \times P', Q \times Q')$ from the oracle and apply either projection. The $\geq$ direction is also easy and we do not need it anyway.)

In our particular case, letting $\mathbf{a}$ and $\mathbf{b}$ be respectively represented by some $A, B \subseteq \N$, we obtain that $\mathbf{c}_{H_0^A \times H_0^B, \, H_1^A \times H_1^B} = \Om$. By \autoref{main-strategy-result},
\[\mathbf{d}_{H_0^A \times H_0^B \times H_1^A \times H_1^B, \, H_1^A \times H_1^B \times H_0^A \times H_0^B} = \mathbf{0}.\]

To simplify this, observe that without any assumptions on $A, B \subseteq \N$, we have in general
\[\mathbf{d}_{H_0^A \times H_0^B, \, H_1^A \times H_1^b} = \mathbf{d}_{H_0^A \times H_0^B \times H_1^A \times H_1^B, \, H_1^A \times H_1^B \times H_0^A \times H_0^B}.\]

Indeed, a possible reduction for $\leq$ takes programs $e_A, e_B$, constructs the oracle program $e'_A$ that runs $e_A$ and exchanges 0 and 1 before returning the result, likewise constructs $e'_B$, queries the oracle on $\langle e_A, e_B, e'_A, e'_B \rangle$ and returns the result. (Again, the inequality $\geq$ is easy and we do not even use it.)

Coming back to our specific situation, we obtain $\mathbf{d}_{H_0^A \times H_0^B, \, H_1^A \times H_1^B} = \mathbf{0}$, and we conclude by \autoref{lemma-T1-meet}.\qedhere
\end{enumerate}
\end{proof}

\section{The internal point of view}\label{section-topos-theory}

As explained in the introduction, Arthur–Nimue–Merlin games were introduced by Kihara to give a more concrete description of Lawvere–Tierney topologies on the effective topos. So far, we have worked entirely at the level of Arthur–Nimue–Merlin games. However, our results exemplify the usefulness of having both points of view. On the one hand, the proofs, especially that of \autoref{main-strategy-result}, involve Arthur–Nimue–Merlin games, and we have not yet found a natural way to explain them purely in terms of Lawvere–Tierney topologies. On the other hand, the idea of the construction, and especially the twist of using $c_{H_0^A, H_1^A}$ instead of $c_{A, A^c}$, comes from the other side, and the purpose of this section is to explain how.

A contrario, we will also show that our main theorem cannot follow merely from general abstract nonsense because although the co-Turing degree construction makes sense in other toposes than the effective topos, it is ill-behaved in general. Even in the effective topos, its extension from Turing degrees to all Arthur–Nimue–Merlin degrees is not injective. Nevertheless, it could be interesting for further research to investigate the construction in toposes related to the effective topos, especially the Kleene–Vesley topos \cite[section 5.4]{vanOostenBook}.

We first recapitulate the connection between Arthur–Nimue–Merlin degrees and Lawvere–Tierney topologies on the effective topos exhibited by Kihara \cite{KiharaLT}. In the interest of exposition, to enable a broader audience able to understand why Arthur–Nimue–Merlin oracles are the most general kind of reasonable oracles, we do not presuppose at this point that the reader understands the effective topos or Lawvere–Tierney topologies in general. We can nevertheless define Lawvere–Tierney topologies on the effective topos as follows. We consider the maps $J : \P(\N) \to \P(\N)$ such that there exist programs
\begin{align*}
\mathtt{ret} &\in \bigcap_{U \subseteq \N} (U \rzimplies J(U))
\\ \mathtt{bind} &\in \bigcap_{U, V \subseteq \N} ((U \rzimplies J(V)) \rzimplies J(U) \rzimplies J(V))
\end{align*}
where we have used the notation
\[(U \rzimplies V) := \{e \in \N \mid \forall n \in U, \varphi_e(n){\downarrow} \in V\}\]
for the set of programs which, given any element of $U$, terminate and return some element of $V$. (The notation is right-associative, like implication: $U \rzimplies V \rzimplies W$ means $U \rzimplies (V \rzimplies W)$.) These maps are endowed with the preorder defined by $J \preceq K$ when there exists a program $\mathtt{convert} \in \bigcap_{U \subseteq \N} (J(U) \rzimplies K(U))$. The partially ordered set of Lawvere–Tierney topologies on the effective topos is obtained by quotienting this preorder by the associated equivalence relation, $J \equiv K$ when $J \preceq K$ and $K \preceq J$.

Informally, we think of $U \subseteq \N$ as the computational problem of “finding an element in $U$”, and $J(U) \subseteq \N$ as ways to solve the problem with the help of some kind of oracle. The existence of $\mathtt{ret}$ and $\mathtt{bind}$ puts two reasonable requirements: if we have already solved the problem, we can turn the solution without help into a solution with help; and if we have a way to solve $U$ with help and a way to use that solution to solve $V$ again with help, then we can compose them to produce a way to solve $V$ with help. The requirements are strongly reminiscent of the definition of a monad in functional programming (and the names $\mathtt{ret}$ and $\mathtt{bind}$ are accordingly borrowed from the Haskell programming language). The connection will be made precise soon: such a map $J$ defines an actual monad when we view it internally in the effective topos.

As a first example, for each Turing oracle $f : \N \to \N$, the following map $J$ represents a Lawvere–Tierney topology on the effective topos: $J(U)$ is the set of programs with oracle $f$ which (run without input) output an element of $U$. The Arthur–Nimue–Merlin games generalize this to obtain all Lawvere–Tierney topologies on the effective topos.

\begin{definition}\label{definition-oracle-topology-correspondence}
For each Arthur–Nimue–Merlin oracle $g : \N \to \P(\P(\N))$, we define $J_g : \P(\N) \to \P(\N)$ as follows: $J_g(U)$ is the set of Arthur strategies for producing an element of $U$ using the oracle $g$ (more formally, the set of Arthur strategies which are part of some winning Arthur–Nimue strategy in the reduction game from the basic oracle $\{U\}$ to $g$).

Conversely, for each map $J : \P(\N) \to \P(\N)$, we define $J^\natural : \N \to \P(\P(\N))$ by $J^\natural(n) := \{U \mid n \in J(U)\}$.
\end{definition}

\begin{theorem}[Kihara {\cite[Corollary 3.6\footnote{Numbered 3.5 in arXiv version}]{KiharaLT}}]\label{theorem-equivalence-anm-ltts}
For all $g$, the map $J_g$ represents a Lawvere–Tierney topology on the effective topos (the programs $\mathtt{ret}$ and $\mathtt{bind}$ as above exist). Furthermore, ${g \mapsto J_g}$ and $J \mapsto J^\natural$ define an equivalence of preordered sets between Arthur–Nimue–Merlin oracles and maps representing Lawvere–Tierney topologies on the effective topos. In other words, they induce inverse isomorphisms between the partially ordered set of Arthur–Nimue–Merlin degrees and the partially ordered set of Lawvere–Tierney topologies on the effective topos.
\end{theorem}

In order to spare the reader the inconvenience of converting between Kihara's and our setup and notations, and also because the proof in \cite{KiharaLT} contains a repairable error, we sketch a proof in \autoref{appendix-equivalence-anm-ltts}.

\begin{remark}\label{remark-inductive-characterization-winning-strategies}
Let $g$ be an Arthur–Nimue–Merlin oracle and let $U \subseteq \N$. The set $J_g(U)$ can be characterized by the following inductive definition (which will turn out useful for computations with the effective topos, in \autoref{ltt-induced-by-inclusion-of-assemblies}):
\begin{itemize}
\item If $s$ is an Arthur strategy that outputs an element of $U$ without consulting the oracle, then $s \in J_g(U)$.
\item If $s$ is an Arthur strategy which starts by making a query $q$ to the oracle, and there exists $R \in g(q)$ such that for all $r \in R$, we have $C(s, r) \in J_g(U)$, then $s \in J_g(U)$. Here, $C(s, r)$ denotes the strategy that gives $r$ to $s$ as a first move and continues simulating it from that point.
\end{itemize}
This is merely a reformulation of the definition: applying countable choice to pick an entire well-founded tree of such $R$ sets gives precisely the data of a Nimue strategy such that the combined Arthur–Nimue strategy is winning.
\end{remark}

From this point on, we assume familiarity with basic topos theory, especially the internal language of an elementary topos, and with the effective topos in particular. A standard reference on the effective topos is van Oosten's monograph on realizability \cite{vanOostenBook}; see also McLarty \cite[chapter 24]{McLarty}, the lecture notes by Phoa \cite{PhoaIntro}, and Hyland's seminal paper \cite{HylandEff}.

Recall that a Lawvere–Tierney topology on an elementary topos is a morphism $j : \Omega \to \Omega$ which is internally a monad on $\Omega$ preordered by implication, viewed as a thin category. Explicitly, the following conditions must be satisfied internally (we omit the equations required for a monad, which are trivially satisfied because the category is thin): $\forall p : \Omega, p \implies j(p)$ (like $\mathtt{ret}$ in Haskell), and $\forall p \, q : \Omega, (p \implies j(q)) \implies (j(p) \implies j(q))$ (like $\mathtt{bind}$ in Haskell). Lawvere–Tierney topologies are preordered by pointwise implication: $j \leq k$ when internally $\forall p : \Omega, j(p) \implies k(p)$. Unfolding the internal language, one recovers the previous definition in the particular case of the effective topos; $\mathbf{0}$ corresponds to the Lawvere–Tierney topology $\id$ while $\Om$ corresponds to $p \mapsto \lnot\lnot p$ and $\top$ corresponds to $p \mapsto \top$.

Returning to the general case, recall that a subobject $B \inject A$ is said to be $j$-dense when internally $\forall x \in A, j(x \in B)$. The following fact is part of the topos theory folklore:

\begin{proposition}\label{proposition-induced-ltt}
For every subobject $B \inject A$ in an elementary topos, there is a least Lawvere–Tierney topology $j$ making the subobject $j$-dense, which we call the Lawvere–Tierney topology induced by $B \inject A$, and it can be internally defined by
\[j(p) := \forall q : \Omega, (p \implies q) \implies (\forall x : A, ((x \in B \implies q) \implies q)) \implies q.\]
\end{proposition}

This formula should be understood as the impredicative encoding of the following inductive definition:

\begin{itemize}
\item If $p$, then $j(p)$.
\item For all $x : A$, if $x \in B$ implies $j(p)$, then $j(p)$.
\end{itemize}

One checks that, internally, $j(p)$ is the least truth value satisfying the above conditions.

\begin{proof}
The closest reference we could find for this exact formula is the work of Ahman and Bauer \cite[Proposition 3.6]{AhmanBauer}, which takes place in a type-theoretic setting that is a bit different from our setup (e.g., the quantification over $B \inject A$ and the quantification implied by “the least” are internal rather than external), but the proof is the same and is in any case a straightforward verification. This is also sketched by Hyland in \cite[Proposition 16.3]{HylandEff} (in the specific case of the effective topos, but this is irrelevant).
\end{proof}

We proceed to unfold this in the effective topos, in a special case that suffices for our needs. We denote an assembly by $(A, E)$ where $A$ is the underlying set and $E : A \to \P(\N)$ is such that each $E(x)$ is inhabited.

\begin{proposition}\label{ltt-induced-by-inclusion-of-assemblies}
Let $(A, E)$ and $(A, F)$ be two assemblies on the same underlying set such that the identity defines a morphism of assemblies $(A, F) \inject (A, E)$. The Lawvere–Tierney topology on the effective topos induced by this subobject is associated (through the correspondence from \autoref{theorem-equivalence-anm-ltts}) to the oracle $g : q \mapsto \{F(x) \mid x \in A, q \in E(x)\}$.
\end{proposition}

In words, when using the oracle $g$, Arthur chooses a realizer $q$, Nimue chooses $x$ such that $q$ realizes $x$ in the assembly $(A, E)$, then Merlin chooses $r$ realizing $x$ in the assembly $(A, F)$.

\begin{proof}
The idea is that externalizing the inductive definition of the induced Lawvere–Tierney topology recovers the inductive characterization of $J_g$ from \autoref{remark-inductive-characterization-winning-strategies}, as observed by Ahman and Bauer \cite[section 5]{AhmanBauer}.

To formally make sense of this, we could work out general facts about the semantics of such inductive definitions (essentially, propositional W-types) in the effective topos. We can also check just this special case by externalizing the formula from \autoref{proposition-induced-ltt}, which gives
\[U \mapsto \bigcap_{V \subseteq \N} \left((U \rzimplies V) \rzimplies \left(\bigcap_{x \in A} (E(x) \rzimplies (F(x) \rzimplies V) \rzimplies V)\right) \rzimplies V\right).\]
and verifying that this is equivalent for $\equiv$ to $J_g$.
\end{proof}

An immediate consequence is a description of $\mathbf{d}_{P, Q}$ and $\mathbf{c}_{P, Q}$ as induced Lawvere–Tierney topologies.

\begin{corollary}
Let $(P, Q)$ be a promise problem. Let $2_{P, Q}$ be the assembly with underlying set $\{0, 1\}$, where elements of $P$ realize 0 and elements of $Q$ realize 1. The identity on the underlying sets define monomorphisms of assemblies $2 \inject 2_{P, Q} \inject \nabla 2$. The subobject $2 \inject 2_{P, Q}$ induces the Lawvere–Tierney topology on the effective topos associated to the Arthur–Nimue–Merlin degree $\mathbf{d}_{P, Q}$, while $2_{P, Q} \inject \nabla 2$ induces the Lawvere–Tierney topology associated to $\mathbf{c}_{P, Q}$.
\end{corollary}

The next step is to define $2 \inject 2_{P, Q} \inject \nabla 2$ in a purely internal way. Working in any topos, let $\Class := \{p : \Omega \mid \lnot\lnot p \implies p\}$ be the object of classical (or $\lnot\lnot$-stable) truth values, and for a Lawvere–Tierney topology $j$, let $\Xi(j) := \{p : \Class \mid j(p \lor \lnot p)\}$ be the object of classical truth values which are $j$-decidable. When $j \leq k$, there is a natural inclusion $\Xi(j) \inject \Xi(k)$. For example, $\Xi(\id)$ is the object of decidable truth values $\{p : \Omega \mid p \lor \lnot p\}$, which is isomorphic to $2$, while $\Xi(\lnot\lnot)$ is simply (isomorphic to) the object of classical truth values $\Class$.

\begin{proposition}
Let $A \subseteq \N$ and let $j$ be the Lawvere–Tierney topology on the effective topos associated to the Turing degree of $A$. The sequence of morphisms $2 \inject 2_{H_0^A, H_1^A} \inject \nabla 2$ in the effective topos is isomorphic to $2 \inject \Xi(j) \inject \Class$ (defined by $2 \cong \Xi(\id) \inject \Xi(j) \inject \Xi(\lnot\lnot) \cong \Class$).
\end{proposition}

\begin{proof}[Proof sketch.]
The isomorphism $\Class \cong \nabla 2$ is standard (it is used without proof in van Oosten's monograph \cite[p. 129]{vanOostenBook}) and the rest is for completely similar reasons, noting that $J_A(\{0\}) = H_0^A$ and $J_A(\{1\}) = H_1^A$ per \autoref{definition-oracle-topology-correspondence}. 
\end{proof}

In summary, the previous discussion establishes the following theorem, which provides some intuition for the construction of $\coT$: “the closer $\Xi(j)$ is to 2, the further it is from $\Class$”.

\begin{theorem}
Let $\mathbf{d}$ be a Turing degree and let $j$ be the associated Lawvere–Tierney topology on the effective topos. The subobject inclusion $2 \inject \Xi(j)$ induces $j$ itself, while $\Xi(j) \inject \Class$ induces the Lawvere–Tierney topology associated to $\coT(\mathbf{d})$.
\end{theorem}

This construction makes sense in an arbitrary elementary topos, but it is ill-behaved in general. For example, consider an elementary topos where the weak law of excluded middle $\forall p : \Omega, \lnot \lnot p \lor \lnot p$ holds, but the law of excluded middle $\forall p : \Omega, p \lor \lnot p$ fails (which exists \cite[Examples D4.6.3]{JohnstoneElephant}). Then the Lawvere–Tierney topologies $\id$ and $\lnot\lnot$ are distinct, but $\Xi(j) \cong \Class$ for all $j$, so both inclusions $2 \inject \Xi(j) \inject \Class$ always induce $\id$.

Staying in the effective topos, we can observe that this construction extends $\coT$ in a natural way to all Arthur–Nimue–Merlin degrees. Explicitly, $\coT([g])$ is $\mathbf{c}_{J_g(\{0\}), J_g(\{1\})}$. However, the resulting map is not injective. For example, consider the basic oracle $\Error_{1/3} := \{\{0, 1\}, \{0, 2\}, \{1, 2\}\}$ (such $\Error_{1/n}$ oracles have been studied in detail by Lee and van Oosten \cite[section 5]{LeeVanOosten}, then by Kihara \cite[section 4.1]{KiharaLT}). The promise problem $(J(\{0\}), J(\{1\}))$ is decidable: given $s \in J(\{0\}) \cup J(\{1\})$, we can run $s$ in parallel on all finite sequences of elements of $\{0, 1, 2\}$ (“explore the tree of reactions to all possible Merlin moves”), and label by $i \in \{0, 1\}$ those nodes such that $s$ returns $i$, and also those nodes with two children labeled by $i$. When the root is labeled by $i$, we know that $s \in J(\{i\})$. Therefore $\coT([\Error_{1/3}]) = \Om$, although $[\Error_{1/3}] > \mathbf{0}$ (by \autoref{remark-when-basic-oracle-is-zero}).

As a concluding remark, we could alternatively define $\Xi(j)$ as the sheafification of 2 for $j$. This and much else will be expanded on in the full version of this paper.

\bibliography{refs}

@article{AhmanBauer,
  title = {Sheaves as oracle computations},
  author = {Ahman, Danel and Bauer, Andrej},
  year = {2026},
  eprint = {2602.22135},
}

@book{CooperBook,
  title = {Computability theory},
  author = {Cooper, S. Barry},
  year = {2004},
  publisher = {Chapman \& Hall/CRC},
  note = {ISBN 9781584882374},
}

@inproceedings{GoldreichPromiseProblems,
  title = {On promise problems: a survey},
  author = {Goldreich, Oded},
  year = {2006},
  booktitle = {Theoretical computer science: essays in memory of Shimon Even},
  editor = {Goldreich, Oded and Rosenberg, Arnold L. and Selman, Alan L.},
  publisher = {Springer Berlin Heidelberg},
  pages = {254--290},
  doi = {10.1007/11685654_12},
}

@inproceedings{HylandEff,
  title = {The effective topos},
  author = {Hyland, J. Martin E.},
  year = {1982},
  booktitle = {The L. E. J. Brouwer centenary symposium: proceedings of the conference held in Noordwijkerhout, 8–13 June, 1981},
  editor = {Troelstra, Anne S. and van Dalen, Dirk},
  publisher = {North-Holland Publishing Company},
  series = {Studies in Logic and the Foundations of Mathematics},
  volume = {110},
  pages = {165--216},
  doi = {10.1016/S0049-237X(09)70129-6},
  url = {https://dpmms.cam.ac.uk/~jmeh1/Research/Pub81-90/hyland-effectivetopos.pdf},
}

@book{JohnstoneElephant,
  title = {Sketches of an elephant: a topos theory compendium},
  author = {Johnstone, Peter T.},
  year = {2002},
  publisher = {Clarendon Press},
  note = {ISBN 9780198534259 (vol. 1) / 9780198515982 (vol. 2)},
}

@article{KiharaLT,
  title = {{Lawvere–Tierney} topologies for computability theorists},
  author = {Kihara, Takayuki},
  year = {2023},
  journal = {Transactions of the American Mathematical Society, series B},
  volume = {10},
  pages = {48--85},
  doi = {10.1090/btran/134},
}

@article{KiharaOracles,
  title = {Rethinking the notion of oracle: a prequel to {Lawvere–Tierney} topologies for computability theorists},
  author = {Kihara, Takayuki},
  year = {2022},
  eprint = {2202.00188},
}

@article{KiharaNg1,
  title = {The {subTuring} degrees},
  author = {Kihara, Takayuki and Ng, Keng M.},
  year = {2025},
  journal = {Canadian Journal of Mathematics},
  pages = {1--26},
  doi = {10.4153/S0008414X25101922},
}

@article{KiharaNg2,
  title = {The game-theoretic {Katětov} order and idealised effective subtoposes},
  author = {Kihara, Takayuki and Ng, Ming},
  year = {2026},
  eprint = {2602.08138},
}

@article{KiharaNg2Bis,
  title = {The gamified {Katětov} order is not linear (in fact, very much not so)},
  author = {Kihara, Takayuki and Ng, Ming},
  year = {2026},
  eprint = {2605.21473},
}

@article{LeeVanOosten,
  title = {Basic subtoposes of the effective topos},
  author = {Lee, Sori and van Oosten, Jaap},
  year = {2013},
  journal = {Annals of Pure and Applied Logic},
  volume = {164},
  issue = {9},
  pages = {866--883},
  doi = {10.1016/j.apal.2013.04.001},
}

@book{McLarty,
  title = {Elementary categories, elementary toposes},
  author = {McLarty, Colin},
  year = {1992},
  publisher = {Oxford University Press},
  doi = {10.1093/oso/9780198533924.001.0001},
  note = {ISBN 9780198533924},
}

@misc{MO507736,
  title = {Answer to “‘{Uniformly}’ minimal pairs of {Turing} degrees” on the {MathOverflow} website},
  author = {Pastebee (pseudonym)},
  year = {2026},
  url = {https://mathoverflow.net/a/507736},
}

@article{Phoa,
  title = {Relative computability in the effective topos},
  author = {Phoa, Wesley},
  year = {1989},
  journal = {Mathematical Proceedings of the Cambridge Philosophical Society},
  volume = {106},
  number = {3},
  pages = {419--422},
  doi = {10.1017/S0305004100068146},
}

@misc{PhoaIntro,
  title = {An introduction to fibrations, topos theory, the effective topos and modest sets},
  author = {Phoa, Wesley},
  url = {https://www.lfcs.inf.ed.ac.uk/reports/92/ECS-LFCS-92-208/},
}

@phdthesis{PittsThesis,
  title = {The theory of triposes},
  author = {Pitts, Andrew M.},
  year = {1981},
  url = {https://www.cl.cam.ac.uk/~amp12/papers/thet/thet.pdf},
}

@book{Soare,
  title = {Turing computability},
  author = {Soare, Robert I.},
  year = {2016},
  publisher = {Springer},
  series = {Theory and Applications of Computability},
  doi = {10.1007/978-3-642-31933-4},
}

@book{vanOostenBook,
  title = {Realizability: an introduction to its categorical side},
  author = {van Oosten, Jaap},
  year = {2008},
  publisher = {Elsevier},
  series = {Studies in Logic and the Foundations of Mathematics},
  volume = {152},
  doi = {10.1016/s0049-237X(08)X8001-2},
  note = {ISBN 9780444515841},
}

@article{vanOostenPittsOperator,
  title = {Realizability with a local operator of {A.} {M.} {Pitts}},
  author = {Jaap van Oosten},
  year = {2014},
  journal = {Theoretical Computer Science},
  volume = {546},
  pages = {237--243},
  doi = {10.1016/j.tcs.2014.03.011},
}

\appendix
\section{Proof of isomorphism between Arthur–Nimue–Merlin degrees and Lawvere–Tierney topologies on the effective topos}\label{appendix-equivalence-anm-ltts}

For convenience of the reader and completeness of the literature, we sketch a proof of \autoref{theorem-equivalence-anm-ltts}, which is due to Kihara. The proof in \cite{KiharaLT} contains a mistake which does not affect the result: in the proof of Theorem 3.5, $\mathtt{b}$ is assumed more uniform than it is. We thank Kihara for promptly acknowledging the issue in private communication and providing helpful explanations. The fix is the “continuation-passing trick” explained below.

This fix in fact already appears in the literature in a hidden way, namely as the external unfolding of Ahman and Bauer's proof that every modality is an oracle modality, which happens internally in a slightly different, type-theoretic setup \cite{AhmanBauer}. This roughly corresponds to the fact that every Lawvere–Tierney topology on an elementary topos is induced by some subobject as in \autoref{proposition-induced-ltt} (the proof is the same as theirs: $j$ is tautologically the least Lawvere–Tierney topology $j'$ such that internally $\forall p : \Omega, j(p) \implies j'(p)$, so is induced by $\{\top\} \inject \{p : \Omega \mid j(p)\}$). However, Ahman and Bauer do not obtain an actual isomorphism because they use a different preorder on oracles, namely morphisms of containers, which correspond to Weihrauch reduction (call-once oracles).
 
\begin{proof}[Proof sketch of \autoref{theorem-equivalence-anm-ltts}]
\begin{itemize}
\item The map $J_g$ represents a Lawvere–Tierney topology: The program $\mathtt{ret}$ exists simply because we can computably convert a natural number $n$ to a program that just returns $n$. For $\mathtt{bind}$, given $f \in (U \rzimplies J_g(V))$ and $e \in J_g(U)$, we can compute (without knowing $U$ and $V$) the Arthur strategy that simulates $e$, yielding a final value $n \in U$, then simulates $\varphi_f(n)$. (A corresponding winning strategy for Nimue simulates a winning strategy associated to $e$, then a winning strategy associated to $\varphi_f(n)$.)

\item The map $g \mapsto J_g$ is monotone: If $f \preceq g$, then $J_f \preceq J_g$ because we can computably convert a program that uses the oracle $f$ to a program using $g$, by replacing each oracle call with a full run of a fixed reduction from $f$ to $g$.

\item The map $J \mapsto J^\natural$ is monotone: If $J \preceq K$, then $J^\natural \preceq K^\natural$ by the following reduction. Arthur has an element $a \in J(U)$ for $U$ unknown to him. He converts it to an element $b \in K(U)$ using a fixed witness of $J \preceq K$, then calls the oracle on $b$. Nimue chooses $U$, which is valid because $b \in K(U)$. Merlin gives Arthur an element of $U$ which Arthur just returns.

\item We have $g \preceq (J_g)^\natural$: Arthur has $n \in \N$ and must find an element of $P$ for some unknown $P \in g(n)$. He constructs the Arthur strategy $s$ which queries the $g$ oracle on $n$ and returns the result. Then he queries the $(J_g)^\natural$ oracle on $s$. We have $s \in J_g(P)$, so Nimue can choose $P$, and Merlin gives Arthur an element of $P$ which Arthur returns.

\item We have $(J_g)^\natural \preceq g$: Arthur has a strategy $s$ which can find an element of some unknown $P$ using the oracle $g$, and simply simulates $s$.

\item We have $J \preceq J_{J^\natural}$: Given $n \in J(U)$ for some unknown $U$, we can compute the Arthur strategy $s$ that calls the ($J^\natural$) oracle on $n$ and returns the result.

\item We have $J_{J^\natural} \preceq J$: This is the delicate point, proved with a “continuation-passing trick”. We witness this using the following program $W$. Fix $\mathtt{ret}$ and $\mathtt{bind}$ for $J$. Our program $W$ receives a strategy $s$ which lets Arthur find an element of some unknown $U$ using the $J^\natural$ oracle. It starts simulating $s$. If $s$ declares a final value $v$ then $W$ returns $\varphi_{\mathtt{ret}}(v)$. If $s$ queries the oracle on a question $q$, then $W$ constructs the program $e$ which takes an input $r \in \N$ and runs $W$ itself on the following strategy (“$C(s,r)$”): give $r$ to $s$ as the first oracle answer and continue simulating it after that point. Then $W$ computes $b := \varphi_{\mathtt{bind}}(e)$ and returns $\varphi_b(q)$. (As usual, the self-reference in the definition of $W$ is justified by Klenee's second recursion theorem.) The correctness is proved by induction on $s \in J_{J^\natural}(U)$ using the inductive characterization of $J_{J^\natural}(U)$ from \autoref{remark-inductive-characterization-winning-strategies}. If $s$ declares the final value $v \in U$, then $W$ returns $\varphi_{\mathtt{ret}}(v)$, which is in $J(U)$ by assumption on $\mathtt{ret}$. Suppose $s$ queries the oracle on $q$ and let $R \in J^\natural(q)$ (i.e., $R$ such that $q \in J(R)$) be a Nimue move as in the inductive characterization. The induction hypothesis ensures that $e \in (R \rzimplies J(U))$ and therefore $b \in (J(R) \rzimplies J(U))$ and thus $\varphi_b(q) \in J(U)$, by assumption on $\mathtt{bind}$.\qedhere
\end{itemize}
\end{proof}

\end{document}